# Chernoff-Savage and Hodges-Lehmann results for Wilks' test of multivariate independence

**Marc Hallin**[1] **and Davy Paindaveine**[*,1]

*Université libre de Bruxelles*

**Abstract:** We extend to rank-based tests of multivariate independence the Chernoff-Savage and Hodges-Lehmann classical univariate results. More precisely, we show that the Taskinen, Kankainen and Oja (2004) normal-score rank test for multivariate independence uniformly dominates – in the Pitman sense – the classical Wilks (1935) test, which establishes the Pitman non-admissibility of the latter, and provide, for any fixed space dimensions $p, q$ of the marginals, the lower bound for the asymptotic relative efficiency, still with respect to Wilks' test, of the Wilcoxon version of the same.

## 1. Introduction

### 1.1. Testing multivariate independence

The problem of testing for independence between two random variables with unspecified densities has been among the very first applications of ranks in statistical inference. Spearman's correlation coefficient was proposed as early as 1904 (Spearman [29]), long before Wilcoxon [33]'s rank sum and signed rank tests for location, and Kendall introduced his rank correlation measure in 1938 (Kendall [18]).

The multivariate version of the same problem, namely, testing independence between two random vectors with unspecified densities, is harder. The first rank-based counterpart to the Gaussian likelihood ratio method of Wilks [34] was developed in Chapter 8 of the monograph by Puri and Sen [25] and, for almost thirty years, has remained the only rank-based solution to the problem. The proposed test, however, is based on componentwise rankings, hence is neither invariant under affine transformations, nor distribution-free – unless of course both random vectors have dimension one, in which case we are back to the traditional context of rank-based tests of bivariate independence.

Of course, the asymptotic null distributions of the Puri and Sen test statistics do not depend on the underlying distributions. This however is a simple consequence of central-limit behavior and studentization – the matrix of the quadratic

*Supported by a Mandat d'Impulsion of the Fonds National de la Recherche Scientique, Communauté française de Belgique.

[1]Institut de Recherche en Statistique, E.C.A.R.E.S., and Département de Mathématique, Université libre de Bruxelles, Campus de la Plaine CP, 210 B-1050 Bruxelles, Belgium, e-mail: mhallin@ulb.ac.be; dpaindav@ulb.ac.be

form involved being a consistent empirical version of its population counterpart. Such property, in the terminology of Randles [27] is only a *weak form* of asymptotic distribution-freeness, the *strong form* being asymptotic equivalence with a genuinely distribution-free random variable.

Alternatives to the Puri and Sen test have been made possible by the recent developments of multivariate concepts of signs and ranks. Gieser and Randles [3] propose a sign test based on Randles [28]'s *interdirections.* Taskinen, Kankainen and Oja [30] also propose a sign test, based on the so-called standardized spatial signs, which is asymptotically equivalent to the Gieser and Randles one under elliptic symmetry assumptions. Multivariate ranks are introduced, along with the signs, in Taskinen, Oja and Randles [32], where multivariate analogs of Spearman's *rho* and Kendall's *tau* (based on Wilcoxon scores) are considered. This is extended, in Taskinen, Kankainen and Oja [31], to arbitrary score functions. All these tests rely on exactly affine-invariant and *strongly* (in Randles' sense) asymptotically distribution-free statistics – this latter property, however, is obtained at the expense of stronger assumptions of elliptical symmetry.

The objective of this paper is to show that Chernoff-Savage and Hodges-Lehmann results can be established for the normal (van der Waerden, say) and Wilcoxon score versions of the Taskinen, Kankainen and Oja [31] tests, establishing their excellent performances with respect to the Gaussian likelihood ratio procedure of Wilks [34].

### *1.2. Chernoff-Savage and Hodges-Lehmann results*

Rank-based inference long has been considered as a somewhat heuristic and heteroclite collection of "quick-and-easy" methods applicable under a broad range of assumptions. The Spearman, Kendall, Wilcoxon, ... tests were not introduced in the context of a general theory, but as isolated tools, to be used when everything else fails. Certainly, their performances were not expected to compare favorably with those of their parametric competitors.

This widespread opinion was partly dispelled by two famous papers – Hodges and Lehmann [15] and Chernoff and Savage [1] – establishing that rank-based methods not only compete very well, but even may outperform their parametric counterparts (quite remarkably, Hotelling and Pabst [17] already show that the asymptotic relative efficiency, under bivariate Gaussian densities, of Spearman with respect to the traditional correlation coefficient, is as high as 0.9119). These papers certainly played a triggering role in the subsequent development of the well-structured theory of rank-based inference associated, mainly, with the name of Hájek, culminating in the monographs by Hájek and Šidák (see Hájek, Šidák and Sen [4]) and Puri and Sen [26].

In their celebrated ".864 theorem", Hodges and Lehmann [15] proved that, in the two-sample location model (but this extends to more general location problems, such as one-sample, *c*-sample, ANOVA, regression, etc.), the Pitman asymptotic relative efficiency (ARE) of Wilcoxon (i.e., linear-score) rank tests with respect to their normal-theory competitor (namely, Student's standard two-sample *t*-test) is never less than .864. In other words, irrespective of the underlying distribution, Wilcoxon tests asymptotically never need more than 13.6% observations more than *t*-tests to achieve the same power.

No less surprising is the Chernoff and Savage [1] result establishing the amazing fact that, in the same class of models, the ARE of van der Waerden (i.e., normal-score) rank tests, still with respect to the corresponding standard Gaussian ones, are always larger than one, and that this minimal value is reached at Gaussian distributions only. It should be stressed that these results deal with the "worst cases": both for the Wilcoxon and the van der Waerden tests, it is possible to show that there is no "best case", in the sense that it is possible to construct a sequence of underlying distributions along which AREs (still with respect to standard Gaussian tests) tend to infinity.

With the extensions of rank-based inference to other models and other inference problems, it has been discovered that these Chernoff-Savage and Hodges-Lehmann results are not just a happy accident specific to *location* problems and *univariate* observations. Hallin [5] showed that the van der Waerden version of the *serial* rank tests proposed by Hallin and Puri [13] also uniformly beats (in the Pitman sense) the corresponding everyday practice parametric Gaussian tests based on autocorrelations. As for the extension of the Hodges-Lehmann [15] result to this time series setup, the lower bound of the Wilcoxon version of those tests, still with respect to the parametric Gaussian tests, was shown to be .856 by Hallin and Tribel [14].

Extensions to (nonserial as well as serial) problems involving *multivariate* observations were recently obtained by Hallin and Paindaveine [7, 8, 10], who showed that their various multivariate van der Waerden rank tests uniformly dominate the corresponding parametric Gaussian procedures in a broad class of problems (culminating in the problem of testing linear restrictions on the parameters of a multivariate general linear model with vector ARMA errors); the Pitman non-admissibility of the associated everyday practice Gaussian tests (one-sample and two-sample Hotelling tests, multivariate $F$-tests, multivariate *portmanteau* and Durbin-Watson tests, etc.) follows. Hallin and Paindaveine [7] (resp., Hallin and Paindaveine [8]) also extended Hodges-Lehmann's result to the multivariate location (resp., serial) setup, providing, for any fixed dimension of the observation space, the lower bound for the AREs of the proposed multivariate linear-score (i.e., Wilcoxon type) rank tests with respect to the parametric Gaussian ones.

Results recently have also been obtained beyond the case of *location* parameters (in a broad sense, including regression and autoregression coefficients): Paindaveine [23] shows that the uniform Pitman dominance of normal-score rank-based procedures over their parametric Gaussian competitors extends to the *shape parameter* of elliptical populations; see Hallin and Paindaveine [11] and Hallin, Oja and Paindaveine [6]. It follows that, for any space dimension $k \geq 2$, the Gaussian maximum likelihood estimator of shape, based on empirical covariances, as well as the corresponding Gaussian likelihood ratio tests, are Pitman-nonadmissible.

In this paper, similar results are established for the problem of testing independence between elliptical random vectors. More precisely, we consider the van der Waerden and Wilcoxon versions of the rank score tests recently proposed by Taskinen, Kankainen and Oja [31]. We first establish, under (multivariate extensions of) the Konijn [19] alternatives, a Chernoff-Savage result, showing that the van der Waerden version of their test uniformly dominates the parametric Gaussian procedure – Wilks [34]'s test – which establishes the Pitman-nonadmissibility of the latter. Similarly, we extend the Hodges-Lehmann [15] ".864 result", providing, for any fixed couple of marginal space dimensions $(p, q)$, the lower bound for the asymptotic relative efficiency of the Taskinen, Kankainen and Oja [31] Wilcoxon test with respect to Wilks.

### *1.3. Outline of the paper*

The paper is organized as follows. In Section 2, we define the notation to be used throughout, describe the problem of testing for independence between elliptically symmetric marginals, and briefly recall the rank score tests developed by Taskinen, Kankainen and Oja [31]. In Section 3, we establish the Pitman non-admissibility of Wilks' test for multivariate independence. The analog of Hodges-Lehmann [15]'s result for the problem under study is derived in Section 4. Section 5 briefly concludes.

## 2. Rank-based tests for multivariate independence

### *2.1. Elliptical symmetry*

Recall that the distribution of a random $k$-vector $\mathbf{X}$ is said to be elliptically symmetric if and only if its probability density function is of the form

$$(2.1) \qquad \underline{f}_{\boldsymbol{\mu},\boldsymbol{\Sigma};f}(\mathbf{x}) := c_{k,f} \, (\det \boldsymbol{\Sigma})^{-1/2} f\left(\left((\mathbf{x}-\boldsymbol{\mu})'\boldsymbol{\Sigma}^{-1}(\mathbf{x}-\boldsymbol{\mu})\right)^{1/2}\right), \quad \mathbf{x} \in \mathbb{R}^k,$$

for some $k$-vector $\boldsymbol{\mu}$ (the centre of the distribution), some symmetric positive definite real $k \times k$ matrix $\boldsymbol{\Sigma} = (\Sigma_{ij})$ with $|\Sigma| = 1$ (the *shape* matrix), and some function $f : \mathbb{R}_0^+ \longrightarrow \mathbb{R}^+$ such that $f > 0$ a.e. and $\mu_{k-1;f} := \int_0^\infty r^{k-1} f(r) \, dr < \infty$ ($c_{k,f}$ is a normalization factor depending on the dimension $k$ and $f$). Denote by $P_k(\boldsymbol{\mu}, \boldsymbol{\Sigma}, f)$ the corresponding distribution.

The shape matrix $\boldsymbol{\Sigma}$ determines the orientation and shape of the equidensity contours associated with $\underline{f}_{\boldsymbol{\mu},\boldsymbol{\Sigma};f}$, which are a family of hyper-ellipsoids centered at $\boldsymbol{\mu}$. The definition we are adopting here involves a determinant-based normalization which, in a sense, can be considered as *canonical* : see Paindaveine [24] and Hallin and Paindaveine [8] for details. Other normalizations are also possible (such as $\Sigma_{11} = 1$ or $\operatorname{tr} \boldsymbol{\Sigma} = k$); the choice of such a normalization does not play any role in the sequel.

The problem of testing for multivariate independence (see Section 2.2 below) is invariant under (block-)affine transformations, and so are all tests considered in this paper. Therefore we can restrict – without loss of generality – to the class of centered spherical distributions, for which $\boldsymbol{\mu}$ and $\boldsymbol{\Sigma}$ do coincide with the origin $\mathbf{0}$ in $\mathbb{R}^k$ and the $k$-dimensional identity matrix $\mathbf{I}_k$, respectively.

Under $P_k(\mathbf{0}, \mathbf{I}_k, f)$, the *radial function* $f$ determines the distribution of $\|\mathbf{X}\|$. More precisely, $\|\mathbf{X}\|$ has probability density function $\tilde{f}_k(r) := (\mu_{k-1;f})^{-1} r^{k-1} f(r) I_{[r>0]}$ ($I_A$ stands for the indicator function of the Borel set $A$); denote by $\tilde{F}_k$ the corresponding distribution function.

The assumption that $\mu_{k-1;f} < \infty$ guarantees that (2.1) is a density. The classical Gaussian procedure for testing multivariate independence – namely, Wilks [34]'s test – requires the underlying distribution to have a finite variance (note that Puri and Sen [25] erroneously require finite moments of order four – an error which is repeated in Muirhead and Waternaux [20]). Consequently, when considering AREs with respect to Wilks' test, we will restrict to radial functions satisfying the stronger condition $\mu_{k+1;f} := \int_0^\infty r^{k+1} f(r) \, dr < \infty$, under which the distribution $P_k(\mathbf{0}, \mathbf{I}_k, f)$ has finite second-order moments. One can associate with each radial function $f$ the *radial function type of* $f$ defined as the class $\{f_a, a > 0\}$, where $f_a(r) := f(ar)$, for all $r > 0$. By affine-invariance, one could restrict to parameter values of the form $(\boldsymbol{\mu}, \boldsymbol{\Sigma}, f) = (\mathbf{0}, \mathbf{I}_k, f_{a_0})$ for which the variances of the associated elliptical distributions are equal to $\mathbf{I}_k$. However, it will be convenient in the sequel to consider all possible radial functions, so that we will only fix $(\boldsymbol{\mu}, \boldsymbol{\Sigma}) = (\mathbf{0}, \mathbf{I}_k)$ and let $f$ range over its radial function type. Some extremely mild smoothness conditions on $f$ – which we throughout will assume to be fulfilled – are required to derive AREs. We refer to Taskinen, Kankainen and Oja [31] for details.

The radial function $f$ is said to be Gaussian if and only if $f = \phi_a$ for some $a > 0$, where $\phi(r) := \exp(-r^2/2)$. Under $\mathrm{P}_k(\mathbf{0}, \mathbf{I}_k, \phi)$, the probability density of $\|\mathbf{X}\|$ is $\tilde{\phi}_k(r) := (2^{(k-2)/2}\Gamma(k/2))^{-1}r^{k-1}\phi(r)I_{[r>0]}$ where $\Gamma(.)$ stands for the Euler gamma function; denote by $\tilde{\Phi}_k$ the corresponding distribution function. Therefore, the distribution of $\|\mathbf{X}\|^2 \stackrel{\mathcal{D}}{=} (\tilde{\Phi}_k^{-1}(U))^2$ (throughout, $U$ stands for a random variable uniformly distributed over $(0,1)$), still under $\mathrm{P}_k(\mathbf{0}, \mathbf{I}_k, \phi)$, is $\chi_k^2$, and the distribution function of $\|\mathbf{X}\|$ is simply $\tilde{\Phi}_k(r) = \Psi_k(r^2)$, where $\Psi_k$ denotes the chi-square distribution function with $k$ degrees of freedom.

### 2.2. Testing for multivariate independence

Consider an i.i.d. sample $(\mathbf{X}'_{11}, \mathbf{X}'_{21})', (\mathbf{X}'_{12}, \mathbf{X}'_{22})', \ldots, (\mathbf{X}'_{1n}, \mathbf{X}'_{2n})'$ of $(p+q)$-dimensional random vectors with the same distribution as $(\mathbf{X}'_1, \mathbf{X}'_2)'$. We are interested in the null hypothesis under which $\mathbf{X}_1$ and $\mathbf{X}_2$, with elliptically symmetric distributions $\mathrm{P}_p(\boldsymbol{\mu}_1, \boldsymbol{\Sigma}_1, f)$ and $\mathrm{P}_q(\boldsymbol{\mu}_2, \boldsymbol{\Sigma}_2, g)$, respectively, are mutually independent. As already mentioned, we can – without loss of generality (since the testing problem under study is invariant under block-diagonal affine transformations) – restrict to centered spherical marginal distributions, that is, to assume $(\boldsymbol{\mu}_1, \boldsymbol{\Sigma}_1) = (\mathbf{0}, \mathbf{I}_p)$ and $(\boldsymbol{\mu}_2, \boldsymbol{\Sigma}_2) = (\mathbf{0}, \mathbf{I}_q)$.

The standard parametric procedure is Wilks [34]'s Gaussian likelihood ratio test ($\phi_{\mathcal{N}}$, say), which rejects the null hypothesis (at asymptotic level $\alpha$) whenever

$$-n \log \frac{|\mathbf{S}|}{|\mathbf{S}_{11}||\mathbf{S}_{22}|} > \chi^2_{pq, 1-\alpha},$$

where

$$\mathbf{S} := \begin{pmatrix} \mathbf{S}_{11} & \mathbf{S}_{12} \\ \mathbf{S}_{21} & \mathbf{S}_{22} \end{pmatrix} := \mathrm{ave}_i \left\{ \left[ \begin{pmatrix} \mathbf{X}_{1i} \\ \mathbf{X}_{2i} \end{pmatrix} - \mathrm{ave}_j \begin{pmatrix} \mathbf{X}_{1j} \\ \mathbf{X}_{2j} \end{pmatrix} \right] \left[ \begin{pmatrix} \mathbf{X}_{1i} \\ \mathbf{X}_{2i} \end{pmatrix} - \mathrm{ave}_j \begin{pmatrix} \mathbf{X}_{1j} \\ \mathbf{X}_{2j} \end{pmatrix} \right]' \right\}$$

stands for the partitioned sample covariance matrix, and where $\chi^2_{pq, 1-\alpha}$ denotes the $\alpha$-upper quantile of the chi-square distribution with $pq$ degrees of freedom.

The Taskinen, Kankainen and Oja [31] rank-based competitors of Wilks' test are defined as follows. Define the standardized subvectors $\mathbf{Z}_{11}, \ldots, \mathbf{Z}_{1n}$ associated with original subvectors $\mathbf{X}_{11}, \ldots, \mathbf{X}_{1n}$ as $\mathbf{Z}_{1i} := \hat{\boldsymbol{\Sigma}}_1^{-1/2}(\mathbf{X}_{1i} - \hat{\boldsymbol{\mu}}_1)$, $i = 1, \ldots, n$, where $\hat{\boldsymbol{\mu}}_1$ and $\hat{\boldsymbol{\Sigma}}_1$ are affine-equivariant root-$n$ consistent estimators of location and scatter, respectively. Consider the corresponding *standardized spatial signs* $\mathbf{U}_{1i} = \mathbf{Z}_{1i}/\|\mathbf{Z}_{1i}\|$ and let $R_{1i}$ denote the rank of $\|\mathbf{Z}_{1i}\|$ among $\|\mathbf{Z}_{11}\|, \ldots, \|\mathbf{Z}_{1n}\|$; note that the scatter matrix estimate $\hat{\boldsymbol{\Sigma}}_1$ can be replaced by any affine-equivariant root-$n$ consistent estimator of the shape matrix, as only the directions and ranks of distances are used in the analysis. The statistics $\mathbf{U}_{2i}$ and $R_{2i}$ are defined in the same way within the sample $\mathbf{X}_{21}, \ldots, \mathbf{X}_{2n}$. Letting $K_1, K_2 : (0,1) \to \mathbb{R}$ be two square-integrable *score functions*, the $(K_1, K_2)$-score version of the rank test statistics for multivariate independence proposed in Taskinen, Kankainen and Oja [31] is

$$(2.2) \quad T_{K_1,K_2} := \frac{npq}{\sigma_{K_1}^2 \sigma_{K_2}^2} \left\| \text{ave}_i \left\{ K_1\left(\frac{R_{1i}}{n+1}\right) K_2\left(\frac{R_{2i}}{n+1}\right) \mathbf{U}_{1i}\,\mathbf{U}_{2i}' \right\} \right\|^2,$$

where $\sigma_K^2 := \int_0^1 K^2(u)\,du$, and $\|\mathbf{A}\|^2 := \text{tr}(\mathbf{A}\mathbf{A}')$ is the squared Frobenius norm of $\mathbf{A}$. Under the null hypothesis of independence (with elliptical marginals), this rank score statistic is asymptotically chi-square with $pq$ degrees of freedom, and the associated test $\phi_{K_1,K_2}$ rejects the null hypothesis (still at asymptotic level $\alpha$) as soon as $T_{K_1,K_2} > \chi^2_{pq,1-\alpha}$.

As we show in the sequel, two particular cases (corresponding to two specific types of score functions) of the above rank score tests exhibit a remarkably good uniform efficiency behavior.

## 3. Pitman non-admissibility of Wilks' test

The van der Waerden (normal-score) version of the rank-based test statistic (2.2) is obtained with the score functions $K_1 = \tilde{\Phi}_p^{-1} = (\Psi_p^{-1})^{1/2}$, $K_2 = \tilde{\Phi}_q^{-1} = (\Psi_q^{-1})^{1/2}$ :

$$T_{vdW} := n \left\| \text{ave}_i \left\{ \tilde{\Phi}_p^{-1}\left(\frac{R_{1i}}{n+1}\right) \tilde{\Phi}_q^{-1}\left(\frac{R_{2i}}{n+1}\right) \mathbf{U}_{1i}\,\mathbf{U}_{2i}' \right\} \right\|^2.$$

In order to compute asymptotic relative efficiencies of the resulting van der Waerden test $\phi_{vdW}$ with respect to Wilks' test $\phi_{\mathcal{N}}$, we must embed the null hypothesis of independence in a model allowing for local (contiguous) alternatives of dependence. Such an embedding is not as obvious as in classical parametric models, where local values of the parameter usually provide the required alternatives. As in Gieser and Randles [3] and Taskinen, Kankainen and Oja [31], we adopt (a multivariate extension of) the local alternatives considered by Konijn [19], of the form

$$(3.1) \quad \begin{pmatrix} \mathbf{X}_{1i} \\ \mathbf{X}_{2i} \end{pmatrix} = \begin{pmatrix} (1-n^{-1/2}\delta)\,\mathbf{I}_p & n^{-1/2}\delta\,\mathbf{M} \\ n^{-1/2}\delta\,\mathbf{M}' & (1-n^{-1/2}\delta)\,\mathbf{I}_q \end{pmatrix} \begin{pmatrix} \mathbf{Y}_{1i} \\ \mathbf{Y}_{2i} \end{pmatrix},$$
$$i = 1,\ldots,n,$$

where $\mathbf{Y}_{1i}$ and $\mathbf{Y}_{2i}$ denote mutually independent random vectors, with elliptic distributions $P_p(\mathbf{0},\mathbf{I}_p,f)$ and $P_q(\mathbf{0},\mathbf{I}_q,g)$, respectively, $\mathbf{M}$ is a non-random matrix with appropriate dimensions, and $\delta \in \mathbb{R}$. In that local model, the null hypothesis of independence reduces to $\mathcal{H}_0 : \delta = 0$. As shown in Taskinen, Kankainen and Oja [31], the asymptotic relative efficiency of the rank score test $\phi_{K_1,K_2}$ with respect to Wilks' test does not depend on $\mathbf{M}$. Their results imply that the asymptotic relative efficiency of the van der Waerden test $\phi_{vdW}$ based on $T_{vdW}$ with respect to Wilks' test, under the sequence of local alternatives in (3.1) is

$$(3.2) \quad \text{ARE}_{q,g}^{p,f}(\phi_{vdW}/\phi_{\mathcal{N}}) = \frac{1}{4p^2q^2}\Big(D_p(\phi,f)C_q(\phi,g) + D_q(\phi,g)C_p(\phi,f)\Big)^2,$$

where, denoting by $\varphi_f(r) := -f'(r)/f(r)$ the optimal location score function associated with the radial function $f$, we let

$$C_k(\phi,f) := \mathrm{E}\Big[\tilde{\Phi}_k^{-1}(U)\,\varphi_f(\tilde{F}_k^{-1}(U))\Big] \quad \text{and} \quad D_k(\phi,f) := \mathrm{E}\Big[\tilde{\Phi}_k^{-1}(U)\,\tilde{F}_k^{-1}(U)\Big].$$

TABLE 1
*AREs of the van der Waerden rank score test $\phi_{vdW}$ for multivariate independence with respect to Wilks' test $\phi_{\mathcal{N}}$, under standard multivariate Student (with 3, 4, 6, and 12 degrees of freedom) and standard multinormal densities, for subvector dimensions $p = 2$ and $q = 1, 2, 3, 4, 6$, and 10, respectively*

| | | $\nu_p$ | | | | |
|---|---|---|---|---|---|---|
| $q$ | $\nu_q$ | 3 | 4 | 6 | 12 | $\infty$ |
| 1 | 3 | 1.378 | 1.295 | 1.266 | 1.281 | 1.339 |
| | 4 | 1.293 | 1.190 | 1.141 | 1.135 | 1.167 |
| | 6 | 1.267 | 1.144 | 1.078 | 1.054 | 1.067 |
| | 12 | 1.285 | 1.141 | 1.058 | 1.019 | 1.016 |
| | $\infty$ | 1.343 | 1.174 | 1.072 | 1.017 | 1.000 |
| 2 | 3 | 1.400 | 1.311 | 1.277 | 1.289 | 1.343 |
| | 4 | 1.311 | 1.204 | 1.152 | 1.144 | 1.174 |
| | 6 | 1.277 | 1.152 | 1.085 | 1.060 | 1.072 |
| | 12 | 1.289 | 1.144 | 1.060 | 1.021 | 1.017 |
| | $\infty$ | 1.343 | 1.174 | 1.072 | 1.017 | 1.000 |
| 3 | 3 | 1.417 | 1.323 | 1.286 | 1.294 | 1.346 |
| | 4 | 1.325 | 1.214 | 1.161 | 1.150 | 1.179 |
| | 6 | 1.286 | 1.159 | 1.091 | 1.065 | 1.076 |
| | 12 | 1.292 | 1.146 | 1.062 | 1.023 | 1.019 |
| | $\infty$ | 1.343 | 1.174 | 1.072 | 1.017 | 1.000 |
| 4 | 3 | 1.430 | 1.332 | 1.292 | 1.298 | 1.348 |
| | 4 | 1.336 | 1.223 | 1.167 | 1.156 | 1.183 |
| | 6 | 1.294 | 1.165 | 1.096 | 1.069 | 1.080 |
| | 12 | 1.295 | 1.149 | 1.064 | 1.024 | 1.020 |
| | $\infty$ | 1.343 | 1.174 | 1.072 | 1.017 | 1.000 |
| 6 | 3 | 1.448 | 1.345 | 1.301 | 1.304 | 1.351 |
| | 4 | 1.353 | 1.236 | 1.177 | 1.163 | 1.189 |
| | 6 | 1.306 | 1.175 | 1.103 | 1.075 | 1.085 |
| | 12 | 1.300 | 1.153 | 1.068 | 1.027 | 1.023 |
| | $\infty$ | 1.343 | 1.174 | 1.072 | 1.017 | 1.000 |
| 10 | 3 | 1.471 | 1.361 | 1.312 | 1.311 | 1.353 |
| | 4 | 1.375 | 1.252 | 1.190 | 1.173 | 1.196 |
| | 6 | 1.323 | 1.188 | 1.114 | 1.084 | 1.092 |
| | 12 | 1.308 | 1.159 | 1.073 | 1.032 | 1.027 |
| | $\infty$ | 1.343 | 1.174 | 1.072 | 1.017 | 1.000 |

Some numerical values of these AREs, under multivariate $t$- and normal distributions, are provided in Table 1. All these values are larger than or equal to 1, and seem to be equal to 1 only if both marginals are multinormal.

Taskinen, Kankainen and Oja [31] point out that it is remarkable that, in the multinormal case, the limiting Pitman efficiency of the van der Waerden score test $\phi_{vdW}$ equals that of Wilks' test. But, it is even more remarkable that, as we shall see, the multinormal case is actually the least favorable one for the van der Waerden procedure. Proposition 1 below indeed states that, as soon as one of the marginals is not Gaussian, the van der Waerden test strictly beats Wilks' test (Table 1 provides an empirical confirmation of this fact). The Pitman non-admissibility of Wilks' test follows.

**Proposition 1.** *For all integers $p, q \geq 1$ and all radial functions $f, g$ such that $\mu_{p+1;f} < \infty$ and $\mu_{q+1;g} < \infty$, we have*

$$\mathrm{ARE}_{q,g}^{p,f}(\phi_{vdW}/\phi_{\mathcal{N}}) \geq 1,$$

*where equality holds iff $f$ and $g$ are Gaussian, with common scale.*

To prove this proposition, we need the following intermediate result; see Paindaveine [22] for an elementary proof, based on an idea developed by Gastwirth and

Wolff [2].

**Lemma 1.** *For all integer $k \geq 1$ and all radial function $f$ satisfying $\mu_{k+1;f} < \infty$, we have $D_k(\phi, f)C_k(\phi, f) \geq k^2$, where equality holds iff $f$ is Gaussian.*

*Proof of Proposition 1.* The proof is based on the decomposition

$$\Big(D_p(\phi, f)C_q(\phi, g) + D_q(\phi, g)C_p(\phi, f)\Big)^2 = A_{q,g}^{p,f} + B_{q,g}^{p,f},$$

where

$$A_{q,g}^{p,f} := 4\, D_p(\phi, f)C_p(\phi, f)D_q(\phi, g)C_q(\phi, g)$$

and

$$B_{q,g}^{p,f} := \Big(D_p(\phi, f)C_q(\phi, g) - D_q(\phi, g)C_p(\phi, f)\Big)^2.$$

It directly follows from Lemma 1 that $A_{q,g}^{p,f} \geq 4p^2q^2$, so that (3.2) entails

$$(3.3) \qquad \mathrm{ARE}_{q,g}^{p,f}(\phi_{vdW}/\phi_{\mathcal{N}}) = \frac{1}{4p^2q^2}(A_{q,g}^{p,f} + B_{q,g}^{p,f}) \geq \frac{1}{4p^2q^2}\, A_{q,g}^{p,f} \geq 1.$$

Let us now show that equality holds iff $f$ and $g$ are Gaussian, with the same variance. Equality in (3.3) indeed requires $A_{q,g}^{p,f} = 4p^2q^2$ and $B_{q,g}^{p,f} = 0$. From Lemma 1, $A_{q,g}^{p,f} = 4p^2q^2$ implies that both $f$ and $g$ are Gaussian ($f = \phi_a$ and $g = \phi_b$, say). Now, since $D_k(\phi, \phi_a) = a^{-1}D_k(\phi, \phi) = a^{-1}k$ and $C_k(\phi, \phi_a) = aC_k(\phi, \phi) = aD_k(\phi, \phi) = ak$ for all $k$, we have $B_{q,\phi_b}^{p,\phi_a} = p^2q^2((b/a) - (a/b))^2$, which is equal to zero iff $a = b$. Consequently, equality holds iff $f = g = \phi_a$, for some $a > 0$. □

## 4. A Hodges-Lehmann result for multivariate independence

We now turn to the Wilcoxon test statistic

$$T_W := \frac{9npq}{(n+1)^4}\, \| \mathrm{ave}_i\, \{R_{1i}\, R_{2i}\, \mathbf{U}_{1i}\, \mathbf{U}_{2i}'\} \|^2,$$

which is associated with the score functions $K_1(u) = K_2(u) = u$, $u \in (0, 1)$. The asymptotic relative efficiency of the corresponding Wilcoxon test $\phi_W$ with respect to Wilks' test $\phi_{\mathcal{N}}$, under the sequence of local alternatives (3.1), is given by

$$\mathrm{ARE}_{q,g}^{p,f}(\phi_W/\phi_{\mathcal{N}}) = \frac{9}{4pq}\Big(D_p(I, f)C_q(I, g) + D_q(I, g)C_p(I, f)\Big)^2,$$

where we let

$$C_k(I, f) := \mathrm{E}\Big[U\, \varphi_f(\tilde{F}_k^{-1}(U))\Big] \quad \text{and} \quad D_k(I, f) := \mathrm{E}\Big[U\, \tilde{F}_k^{-1}(U)\Big];$$

see Taskinen, Kankainen and Oja [31]. Some numerical values of these AREs, under multivariate $t$- and normal distributions, are provided in Table 2. The uniformly good asymptotic efficiency of the Wilcoxon test in Table 2 is not just an empirical finding, as shown by the following result which provides the lower bound of these AREs for any fixed values of the dimensions $p, q$ of the marginals (some numerical values of this lower bound are presented in Table 3).

TABLE 2
*AREs of the Wilcoxon rank score test $\phi_W$ for multivariate independence with respect to Wilks' test $\phi_\mathcal{N}$, under standard multivariate Student (with 3, 4, 6, and 12 degrees of freedom) and standard multinormal densities, for subvector dimensions $p = 2$ and $q = 1, 2, 3, 4, 6$ and 10, respectively*

| | | $\nu_p$ | | | | |
|---|---|---|---|---|---|---|
| $q$ | $\nu_q$ | 3 | 4 | 6 | 12 | $\infty$ |
| 1 | 3 | 1.305 | 1.227 | 1.193 | 1.193 | 1.222 |
| | 4 | 1.239 | 1.147 | 1.099 | 1.085 | 1.098 |
| | 6 | 1.208 | 1.104 | 1.044 | 1.018 | 1.018 |
| | 12 | 1.204 | 1.086 | 1.015 | 0.978 | 0.969 |
| | $\infty$ | 1.219 | 1.087 | 1.006 | 0.959 | 0.940 |
| 2 | 3 | 1.305 | 1.237 | 1.211 | 1.219 | 1.257 |
| | 4 | 1.237 | 1.152 | 1.111 | 1.102 | 1.121 |
| | 6 | 1.211 | 1.111 | 1.055 | 1.033 | 1.037 |
| | 12 | 1.219 | 1.102 | 1.033 | 0.997 | 0.989 |
| | $\infty$ | 1.257 | 1.121 | 1.037 | 0.989 | 0.970 |
| 3 | 3 | 1.274 | 1.213 | 1.193 | 1.206 | 1.248 |
| | 4 | 1.203 | 1.125 | 1.089 | 1.084 | 1.106 |
| | 6 | 1.179 | 1.084 | 1.032 | 1.013 | 1.020 |
| | 12 | 1.192 | 1.079 | 1.013 | 0.980 | 0.973 |
| | $\infty$ | 1.245 | 1.110 | 1.027 | 0.979 | 0.960 |
| 4 | 3 | 1.248 | 1.192 | 1.175 | 1.191 | 1.235 |
| | 4 | 1.174 | 1.101 | 1.068 | 1.066 | 1.090 |
| | 6 | 1.149 | 1.059 | 1.011 | 0.993 | 1.002 |
| | 12 | 1.165 | 1.056 | 0.992 | 0.961 | 0.955 |
| | $\infty$ | 1.228 | 1.095 | 1.013 | 0.966 | 0.947 |
| 6 | 3 | 1.211 | 1.161 | 1.150 | 1.168 | 1.215 |
| | 4 | 1.134 | 1.067 | 1.038 | 1.039 | 1.066 |
| | 6 | 1.105 | 1.022 | 0.978 | 0.963 | 0.974 |
| | 12 | 1.122 | 1.019 | 0.959 | 0.930 | 0.927 |
| | $\infty$ | 1.198 | 1.068 | 0.988 | 0.943 | 0.924 |
| 10 | 3 | 1.173 | 1.129 | 1.121 | 1.144 | 1.193 |
| | 4 | 1.090 | 1.029 | 1.005 | 1.009 | 1.038 |
| | 6 | 1.056 | 0.979 | 0.940 | 0.929 | 0.941 |
| | 12 | 1.069 | 0.973 | 0.918 | 0.892 | 0.891 |
| | $\infty$ | 1.158 | 1.033 | 0.955 | 0.911 | 0.893 |

TABLE 3
*Some numerical values, for various values of the dimensions $p$ and $q$ of the subvectors, of the Hodges-Lehmann lower bound for the asymptotic relative efficiency of the Wilcoxon rank score test $\phi_W$ for multivariate independence with respect to Wilks' test $\phi_\mathcal{N}$*

| $p/q$ | 1 | 2 | 3 | 4 | 6 | 10 | $\infty$ |
|---|---|---|---|---|---|---|---|
| 1 | 0.856 | 0.884 | 0.867 | 0.850 | 0.826 | 0.797 | 0.694 |
| 2 | | 0.913 | 0.895 | 0.878 | 0.853 | 0.823 | 0.717 |
| 3 | | | 0.878 | 0.861 | 0.836 | 0.807 | 0.703 |
| 4 | | | | 0.845 | 0.820 | 0.792 | 0.689 |
| 6 | | | | | 0.797 | 0.769 | 0.669 |
| 10 | | | | | | 0.742 | 0.646 |
| $\infty$ | | | | | | | 0.563 |

**Proposition 2.** *Let $p, q \geq 1$ be two integers. Then, letting*

$$c_k := \inf\left\{ x > 0 \,\bigg|\, \left(\sqrt{x}\, J_{\sqrt{2k-1}/2}(x)\right)' = 0 \right\}, \quad k \in \mathbb{N}_0,$$

*where $J_r$ denotes the first-kind Bessel function of order $r$, the lower bound for the asymptotic relative efficiency of $\phi_W$ with respect to $\phi_\mathcal{N}$, for fixed subvector*

*dimensions $p, q$, is*

$$\inf_{f,g} \mathrm{ARE}_{q,g}^{p,f}(\phi_W/\phi_{\mathcal{N}}) = \frac{9}{2^{10} pq c_p^2 c_q^2} \left(2 c_p^2 + p - 1\right)^2 \left(2 c_q^2 + q - 1\right)^2,$$

*where the infimum is taken over the collection of radial functions $f, g$ for which $\mu_{p+1;f} < \infty$ and $\mu_{q+1;g} < \infty$. The infimum is reached at*

$$(f, g) \in \left\{ \left(h_{p,\sigma}(r), h_{q,\sigma}(r)\right) := \left(h_{p,1}(\sigma r), h_{q,1}(\sigma r)\right), \sigma > 0 \right\},$$

*where $h_{k,1}$ denotes "the" ($h_{k,1}$ is defined up to a positive scalar factor) radial function associated with the radial cumulative distribution function*

$$H_{k,1}(r) := \frac{\sqrt{r}\, J_{\sqrt{2k-1}/2}(r)}{\sqrt{c_k}\, J_{\sqrt{2k-1}/2}(c_k)}\, I_{[0 < r \leq c_k]} + I_{[r > c_k]}.$$

To prove this proposition, we need the following lemma, which is established in the proof of Proposition 7 of Hallin and Paindaveine [8].

**Lemma 2.** *Let $k \geq 1$ be a fixed integer. Then,*

$$\inf_f \left\{ D_k(I, f) C_k(I, f) \right\} = \frac{1}{2^5 c_k^2} \left(2 c_k^2 + k - 1\right)^2,$$

*where the infimum is taken over the collection of radial functions $f$ for which $\mu_{k+1;f} < \infty$, and the infimum is reached at the radial functions $f \in \{h_{k,\sigma}(r), \sigma > 0\}$. Moreover, letting $\omega_k := (2 c_k^2 + k - 1)/(8 c_k)$, we have $D_k(I, h_{k,\omega_k}) = 1$.*

Since we have $D_k(I, f_a) = a^{-1} D_k(I, f)$ and $C_k(I, f_a) = a C_k(I, f)$ for all $k$, the quantity $D_k(I, f_a) C_k(I, f_a)$ does not depend on $a$. This allows for identifying a particular member $h_{k,\sigma_k}$ of the radial function type $\{h_{k,\sigma}(r), \sigma > 0\}$ such that $D_k(I, h_{k,\sigma_k}) = 1$. According to Lemma 2, $\sigma_k = \omega_k$.

*Proof of Proposition 2.* Proceeding as in the proof of Proposition 1, we consider the decomposition

$$\mathrm{ARE}_{q,g}^{p,f}(\phi_W/\phi_{\mathcal{N}}) = A_{q,g}^{p,f} + B_{q,g}^{p,f},$$

where

$$A_{q,g}^{p,f} := \frac{9}{pq} D_p(I, f) C_p(I, f) D_q(I, g) C_q(I, g),$$

$$B_{q,g}^{p,f} := \frac{9}{4pq} \left(D_p(I, f) C_q(I, g) - D_q(I, g) C_p(I, f)\right)^2.$$

Lemma 2 directly yields that, for all couple $(f, g)$ of radial functions,

(4.1) $\quad \mathrm{ARE}_{q,g}^{p,f}(\phi_W/\phi_{\mathcal{N}}) \geq A_{q,g}^{p,f} \geq \dfrac{9}{2^{10} pq c_p^2 c_q^2} \left(2 c_p^2 + p - 1\right)^2 \left(2 c_q^2 + q - 1\right)^2.$

We now show that the right-hand side in (4.1) actually coincides with the infimum, by determining the (non-empty) collection of couples $(f, g)$ achieving the bound in (4.1). At the couple $(f, g)$, the bound is achieved iff

(4.2) $\quad A_{q,g}^{p,f} \;=\; \dfrac{9}{2^{10} pq c_p^2 c_q^2} \left(2 c_p^2 + p - 1\right)^2 \left(2 c_q^2 + q - 1\right)^2, \quad \text{and}$

$\qquad\qquad B_{q,g}^{p,f} \;=\; 0.$

Lemma 2 shows that (4.2) holds iff

$$(f,g) \in \left\{ \left(h_{p,a}(r), h_{q,b}(r)\right) := \left(h_{p,1}(ar), h_{q,1}(br)\right), a,b > 0 \right\}.$$

Now, using the fact that $D_k(I, h_{k,1}) = \omega_k$ for all $k$, we have

$$\begin{aligned}
B^{p,h_{p,a}}_{q,h_{q,b}} &= \frac{9}{4pq}\left(\frac{b}{a}D_p(I,h_{p,1})C_q(I,h_{q,1}) - \frac{a}{b}D_q(I,h_{q,1})C_p(I,h_{p,1})\right)^2 \\
&= \frac{9}{4pq}\left(\frac{\omega_p b}{a}C_q(I,h_{q,1}) - \frac{\omega_q a}{b}C_p(I,h_{p,1})\right)^2 \\
&= \frac{9}{4pq}\left(\frac{\omega_p b}{\omega_q a}D_q(I,h_{q,1})C_q(I,h_{q,1}) - \frac{\omega_q a}{\omega_p b}D_p(I,h_{p,1})C_p(I,h_{p,1})\right)^2 \\
&= 0
\end{aligned}$$

iff

$$\left(\frac{\omega_q a}{\omega_p b}\right)^2 = \frac{D_q(I,h_{q,1})C_q(I,h_{q,1})}{D_p(I,h_{p,1})C_p(I,h_{p,1})} = \frac{c_p^2(2c_q^2+q-1)^2}{c_q^2(2c_p^2+p-1)^2} = \left(\frac{\omega_q}{\omega_p}\right)^2,$$

that is, iff $a = b$. □

Somewhat surprisingly, the main results of this paper are scale-dependent. Essentially, the bounds in Propositions 1 and 2 are achieved at densities $f$ and $g$ with common scale. This at first sight is puzzling. Indeed, all tests involved are (block-)affine-invariant. This dependence on scale however is entirely due to the form of the local Konijn alternatives (3.1) considered throughout. Appropriate rescaling of **M** would take care of it. For the sake of coherence with the literature (see, e.g., Gieser and Randles [3], Taskinen, Kankainen and Oja [30, 31] and Taskinen, Oja and Randles [32]), however, we sticked to the traditional definition of Konijn alternatives.

## 5. Final comments

Multivariate signed rank tests, based on affine-invariant concepts of signs and ranks – Randles' interdirections (Randles [28]); (pseudo)-Mahalanobis ranks (Hallin and Paindaveine [7, 8, 9, 10, 11]; Hallin, Oja and Paindaveine [6]); Tyler's angles (Hallin and Paindaveine [9]); Oja-Paindaveine ranks (Oja and Paindaveine [21]) – recently have attracted renewed attention to rank-based methods in the context of multivariate analysis, a domain where the 1971 monograph by Puri and Sen has remained the main reference for more than thirty years.

The specific problem of testing independence between elliptical random vectors has been investigated in a series of papers by Taskinen, Kankainen and Oja [30, 31] and Taskinen, Oja and Randles [32]. Simulations show that the performances of the procedures proposed in these papers are quite good. This empirical finding is confirmed here by an investigation of their asymptotic relative efficiencies with respect to their traditional Gaussian counterpart, Wilks' test of independence, based on classical correlations. More specifically, we obtain, for the Wilcoxon and van der Waerden versions of the Taskinen-Kankainen-Oja tests, analogues of the classical Chernoff-Savage [1] and Hodges-Lehmann [15] results. The Chernoff-Savage result of Proposition 1 in particular establishes the non-admissibility of Wilks' procedure, which is uniformly dominated, in the Pitman sense, by the van der Waerden version of Taskinen, Kankainen and Oja [31].

**Acknowledgments.** Part of this work was completed as Marc Hallin was visiting the Department of Operations Research and Financial Engineering and Bendheim Center at Princeton University. Their financial support and congenial research environment is gratefully acknowledged.

## References

[1] CHERNOFF, H. AND SAVAGE, I. R. (1958). Asymptotic normality and efficiency of certain nonparametric tests. *Ann. Math. Statist.* **29** 972–994. MR0100322

[2] GASTWIRTH, J. L. AND WOLFF, S. S. (1968). An elementary method for obtaining lower bounds on the asymptotic power of rank tests. *Ann. Math. Statist.* **39** 2128–2130. MR0233472

[3] GIESER, P. W. AND RANDLES, R. H. (1997). A nonparametric test of independence between two vectors. *J. Amer. Statist. Assoc.* **92** 561–567. MR1467849

[4] HÁJEK, I., ŠIDÁK, Z. AND SEN, P. K. (1999). *Theory of Rank Tests*, 2nd ed. Academic Press, New York. MR1680991

[5] HALLIN, M. (1994). On the Pitman-nonadmissibility of correlogram-based methods. *J. Time Ser. Anal.* **15** 607–612. MR1312324

[6] HALLIN, M., OJA, H. AND PAINDAVEINE, D. (2006). Semiparametrically efficient rank-based inference for shape: II. Optimal R-estimation of shape. *Ann. Statist.* **34** 2757–2789. MR2329466

[7] HALLIN, M. AND PAINDAVEINE, D. (2002a). Optimal tests for multivariate location based on interdirections and pseudo-Mahalanobis ranks. *Ann. Statist.* **30** 1103–1133. MR1926170

[8] HALLIN, M., AND PAINDAVEINE, D. (2002b). Optimal procedures based on interdirections and pseudo-Mahalanobis ranks for testing multivariate elliptic white noise against ARMA dependence. *Bernoulli* **8** 787–816. MR1963662

[9] HALLIN, M. AND PAINDAVEINE, D. (2002c). Multivariate signed ranks: Randles' interdirections or Tyler's angles? In *Statistical Data Analysis Based on the $L_1$ Norm and Related Procedures* (Y. Dodge, ed.) 271–282. Birkhäuser, Basel. MR2001322

[10] HALLIN, M. AND PAINDAVEINE, D. (2005). Affine invariant aligned rank tests for the multivariate general linear model with ARMA errors. *J. Multivariate Anal.* **93** 122–163. MR2119768

[11] HALLIN, M. AND PAINDAVEINE, D. (2006a). Semiparametrically efficient rank-based inference for shape. I. Optimal rank-based tests for sphericity. *Ann. Statist.* **34** 2707–2756. MR2329465

[12] HALLIN, M. AND PAINDAVEINE, D. (2006b). Parametric and semiparametric inference for shape: The role of the scale functional. *Statist. Decisions* **24** 1001–1023. MR2305111

[13] HALLIN, M. AND PURI, M. L. (1994). Aligned rank tests for linear models with autocorrelated error terms. *J. Multivariate Anal.* **50** 175–237. MR1293044

[14] HALLIN, M. AND TRIBEL, O. (2000). The efficiency of some nonparametric rank-based competitors to correlogram methods. In *Game Theory, Optimal Stopping, Probability, and Statistics* (F. T. Bruss and L. Le Cam, eds.). *Papers in Honor of T.S. Ferguson on the Occasion of his 70th Birthday. I.M.S. Lecture Notes-Monograph Series* 249–262. MR1833863

[15] HODGES, J. L. AND LEHMANN, E. L. (1956). The efficiency of some nonparametric competitors of the *t*-test. *Ann. Math. Statist.* **27** 324–335. MR0079383

[16] HOTELLING, H. (1931). The generalization of Student's ratio. *Ann. Math. Statist.* **2** 360–378.
[17] HOTELLING, H. AND PABST, M. R. (1936). Rank correlation and tests of significance involving no assumption of normality. *Ann. Math. Statist.* **7** 29–43.
[18] KENDALL, M. G. (1938). A new measure of rank correlation. *Biometrika* **30** 81–93.
[19] KONIJN, H. S. (1956). On the power of certain tests for independence in bivariate populations. *Ann. Math. Statist.* **27** 300–323. MR0079384
[20] MUIRHEAD, R. J. AND WATERNAUX, C. M. (1980). Asymptotic distributions in canonical correlation analysis and other multivariate procedures for nonnormal populations. *Biometrika* **67** 31–43. MR0570502
[21] OJA, H. AND PAINDAVEINE, D. (2005). Optimal signed-rank tests based on hyperplanes. *J. Statist. Plann. Inference* **135** 300–323. MR2200471
[22] PAINDAVEINE, D. (2004). A unified and elementary proof of serial and nonserial, univariate and multivariate, Chernoff-Savage results. *Statist. Methodology* **1** 81–91. MR2160622
[23] PAINDAVEINE, D. (2006). A Chernoff-Savage result for shape. On the nonadmissibility of pseudo-Gaussian methods. *J. Multivariate Anal.* **97** 2206–2220. MR2301635
[24] PAINDAVEINE, D. (2008). A canonical definition of shape. *Statist. Probab. Lett.* To appear.
[25] PURI, M. L. AND SEN, P. K. (1971). *Nonparametric Methods in Multivariate Analysis.* Wiley, New York. MR0298844
[26] PURI, M. L. AND SEN, P. K. (1985). *Nonparametric Methods in General Linear Models.* Wiley, New York. MR0794309
[27] RANDLES, R. H. (1984). On tests applied to residuals. *J. Amer. Statist. Assoc.* **79** 349–354. MR0755091
[28] RANDLES, R. H. (1989). A distribution-free multivariate test based on interdirections. *J. Amer. Statist. Assoc.* **84** 1045–1050. MR1134492
[29] SPEARMAN, C. (1904). The proof and measurement of association between two things. *Amer. J. Psychology* **15** 72–101.
[30] TASKINEN, S., KANKAINEN, A. AND OJA, H. (2003). Sign test of independence between two random vectors. *Statist. Probab. Lett.* **62** 9–21. MR1965367
[31] TASKINEN, S., KANKAINEN, A. AND OJA, H. (2004). Rank scores tests of multivariate independence. In *Theory and Applications of Recent Robust Methods* (M. Hubert, G. Pison, A. Struyf and S. Van Aelst, eds.) 329–342. Birkhäuser, Basel. MR2088309
[32] TASKINEN, S., OJA, H. AND RANDLES, R. (2005). Multivariate nonparametric tests of independence. *J. Amer. Statist. Assoc.* **100** 916–925. MR2201019
[33] WILCOXON, F. (1945). Individual comparisons by ranking methods. *Biometrics Bulletin* **1** 80–83.
[34] WILKS, S. S. (1935). On the independence of $k$ sets of normally distributed statistical variables. *Econometrica* **3** 309–326.